# Capacitated Assortment Optimization with Pricing under the Paired Combinatorial Logit Model

Daihan Zhang[1], Zhenghe Zhong[1], Chuning Gao[1] , Rui Chen[2],
[1]Sparkzone Institute, Beijing, 10084, China
[2]Department of Industrial Engineering, Tsinghua University, Beijing 100084, China
E-mail Address: jackdaihan.zhang@outlook.com (Daihan Zhang)
484245118@qq.com (Zhenghe Zhong)
13718121330@163.com (Chuning Gao)
chenr15@mails.tsinghua.edu.cn (Rui Chen)

*Abstract* - **In this paper, we investigate the capacitated assortment optimization problem with pricing under the paired combinatorial logit model, whose goal is to identify the revenue-maximizing subset of products as well as their selling prices subject to a known capacity limit. We model customers' purchase behavior using the paired combinatorial logit model, which allows for covariance among any pair of products. We formulate this problem as a non-linear mixed integer program. Then, we propose a two-step approach to obtain the optimal solution based on solving a mixed integer program and Lambert-W function. To further improve its performance, we design a greedy heuristic algorithm and a greedy randomized adaptive search procedure to obtain high-quality solutions so as to balance the tradeoff between accuracy and computational efficiency. A series of numerical experiments are conducted to gauge the efficiency and quality of our proposed approaches.**

*Keywords* - assortment optimization, price optimization, paired combinatorial logit, revenue management

## I. INTRODUCTION

In the field of revenue management, retailers who sell several categories of products, usually have flexibility on determining the subset of products to offer as well as their selling prices to some extent. Therefore, given a known capacity constraint to limit the total shelf space consumption or the total costs of introducing products to the store, these retailers face the capacitated assortment optimization problem with pricing, referred to as the CAOPP, to maximize the expected revenue obtained from customers. Take a home appliance store as an example, the retailer has to choose a revenue-maximizing subset of televisions, refrigerators, and other products with their corresponding prices of different sizes to display as the total shelf space is limited.

To tackle this problem, researchers mainly introduce variants of discrete choice models, especially the multinomial logit and the nested logit models [see, e.g., 1-3], to capture customers' purchase behavior as well as the substitution among products. Under the multinomial logit model, several efficient algorithms are proposed to solve the corresponding CAOPP with different side constraints [1-2, 4]. Reference [5] characterizes the equilibrium in a competitive environment with multiple retailers. Although analytically convenient, the multinomial logit model suffers from the independence of irrelevant alternative property, which is incapable of handling with the substitution among products [6]. Researchers, therefore, resort to the nested logit model developed by [7]. A sampling of these research on the CAOPP includes [3] and [8-10]. We refer interested readers to [11] for study with rank-based nonparametric model.

In this paper, we investigate the CAOPP under the paired combinatorial logit (PCL) model, which has been proven to be an effective model in capturing the decision making process. There are several benefits of the PCL model: (a) The PCL model is consistent with the random utility maximization principle; (b) The PCL model allows for covariance among any pair of products, which leads to a more accurate representation of choice setting without specifying a structural sequence; and (c) Most, if not all, of extant empirical studies point out that the PCL model outperforms the multinomial logit and the nested logit models in predicting users' route choice [see, e.g., 12-15]. However, there are very scarce applications involving the PCL model in the field of revenue management. Reference [16] is the first to study the pricing problem under the PCL model. The authors show that the uniqueness of the optimal prices can be achieved under certain conditions based on the concept of P-matrix and develop the corresponding solution algorithms. Recently, Reference [17] shows that the assortment optimization problem under the PCL model is NP-hard even when there are no capacity constraints and develop a general approximation framework based on iterative variable fixing and coupled randomized rounding. Later, Reference [18] proposes another approximate algorithm based on approximations to the knapsack problem when there is a known capacity constraint. As far as the authors have concerned, no literature has ever contributed to the CAOPP under the PCL model, and our research fills this blank.

We first formulate this problem as a mixed integer program involving a non-linear objective function. Then, we propose a two-step approach to obtain the optimal solution based on the observation that the optimal prices for all offered products are identical. In the first step, we solve a non-linear auxiliary problem and provide its equivalent integer program, whose linear programming relaxation provide a numerically tight upper bound. In the second step, we obtain the optimal prices for all products based on Lambert-W function and the output in the first step. To further improve its empirical performance, we apply a greedy heuristic algorithm and a randomized adaptive search procedure (GRASP) to obtain high-quality solutions. Numerical experiments indicate that our greedy heuristic algorithm is quite efficient. The optimality gap of the solutions obtained by the heuristic algorithm is less than 0.18% on average and 0.46% at the worst case over 300 randomly generated problem instances.

We summarize the contributions of this research as follows: (a) To the best of authors' knowledge, we are the first to investigate the capacitated assortment optimization problem with pricing under the paired combinatorial logit model; (b) We obtain the optimal assortment as well as the optimal prices through solving a mixed integer program and a Lambert-W function; (c) To further improve its performance, we proceed to develop a greedy heuristic algorithm and a GRASP to obtain high-quality solutions within reasonable time consumed; and (d) We demonstrate the efficiency and quality of our approaches through a series of numerical experiments.

The remainder of this paper is organized as follows. In Section II, we formulate the CAOPP under the PCL model as a mixed integer program. In Section III, we develop a two-step approach based on an integer program to obtain the optimal solution, the upper bound of the original optimization problem, and to develop a greedy heuristic algorithm as well as a GRASP to obtain high-quality solutions quickly. In Section IV, we conduct a series of numerical experiments to gauge the performance of our approaches. Finally, we conclude and outline potential future research directions in Section V.

## II. MODELING FRAMEWORK

In this section, we formulate our joint capacitated assortment and price optimization problem under the PCL model. Consider a retailer, who has enough flexibility on determining the subset of products as well as their selling prices out with limited display space for offered products. Without loss of generality, we assume that the set of all possible products to offer is indexed by $N = \{1, 2, \cdots, n\}$. For each product $i \in N$, let $p_i$ be its selling price and $w_i$ be its weight, or required display space. Moreover, we use $x_i = 1$ to indicate that product $i$ is offered; otherwise, $x_i = 0$. Given a capacity limit, denoted by $C \in \mathbb{R}^+$, to limit the total space consumption of the offered products, the retailer wishes to find the revenue-maximizing subset of products as well as their selling prices.

In this paper, we assume that customers' choice behavior is captured by the PCL model, which can accurately capture the substitution effects among products [16-18]. Throughout this paper, the deterministic utility of product $i$ is captured through $\alpha_i - \beta p_i$, where $\alpha_i$ is a known constant to measure the quality of this product, and $\beta > 0$ is the price-sensitivity parameter to capture the utility variation when the selling price changes. Under the PCL model, products are partitioned into several nests which contain exact two products, and we use $\langle i, j \rangle (i < j)$ to denote the nest with only products $i$ and $j$. Let $\mathbf{p} = (p_1, p_2, \cdots, p_n)$ and $\mathbf{x} = (x_1, x_2, \cdots x_n)$. We use $\gamma_{ij} \in (0,1]$ to denote the dissimilarity parameter associated with nest $\langle i, j \rangle$, the choice probability of purchasing product $i$ can be calculated through $q_i^{ij} = v_i^{1/\gamma_{ij}}(p_i)x_i / V_{ij}(\mathbf{p}, \mathbf{x})$, where $v_i = \exp(\alpha_i - \beta p_i)$ is the preference weight of product $i$ and $V_{ij}(\mathbf{p}, \mathbf{x}) = v_i^{1/\gamma_{ij}}(p_i)x_i + v_j^{1/\gamma_{ij}}(p_j)x_j$ is the total preference weights of nest $\langle i, j \rangle$. Moreover, the convention that $q_i^{ij} = 0$ is adopted if $x_i = x_j = 0$. Then, the expected revenue that we obtain from this customer is $R_{ij}(\mathbf{p}, \mathbf{x}) = p_i q_i^{ij} + p_j q_j^{ij}$. Furthermore, the choice probability $q^{ij}$ that a customer will make a purchase under nest $\langle i, j \rangle$ is $q^{ij} = V_{ij}^{\gamma_{ij}}(\mathbf{p}, \mathbf{x}) / (1 + \sum_{i=1}^{n-1} \sum_{j=i+1}^{n} V_{ij}^{\gamma_{ij}}(\mathbf{p}, \mathbf{x}))$. Finally, the choice probability of product $i$, denoted by $q_i(\mathbf{p}, \mathbf{x})$ with given $\mathbf{x}$ and $\mathbf{p}$, is $q_i(\mathbf{p}, \mathbf{x}) = \sum_i^n \sum_j^n q^{ij} \cdot q_i^{ij}$, and the no-purchase option, which means the customer leaves without purchasing anything, is $q_0 = 1 - \sum_{i=1}^{n-1} \sum_{j=i+1}^{n} q^{ij}$. As a result, the CAOPP under the PCL model can be formulated as the following non-linear mixed integer program:

$$z^* = \max R(\mathbf{p},\mathbf{x}) = \sum_{i=1}^{n} p_i \cdot q_i = \sum_{i=1}^{n-1} \sum_{j=i+1}^{n} p_i \cdot q^{ij} \cdot q_i^{ij}$$

$$= \sum_{i=1}^{n-1} \sum_{j=i+1}^{n} \frac{V_{ij}^{\gamma_{ij}}(\mathbf{p},\mathbf{x}) R_{ij}(\mathbf{p},\mathbf{x})}{1 + \sum_{i=1}^{n-1} \sum_{j=i+1}^{n} V_{ij}^{\gamma_{ij}}(\mathbf{p},\mathbf{x})} \quad (1)$$

s.t. $\sum_{i=1}^{n} w_i x_i \leq C,$ (2)

$p_i \geq 0, \ x_i \in \{0,1\}, \ \forall i \in N.$ (3)

The objective function maximizes the expected revenue over all products. Constraint (2) ensures that the total weight consumption of offered products will not exceed the capacity limit. Constraints (3) restrict the decision variables.

## III. SOLUTION APPROACH

In this section, we present our solution approaches to obtain the optimal solution as well as high-quality solutions to the problem defined in (1)-(3).

The starting point of our solution approaches is a result from the observation that with any given feasible assortment $\mathbf{x}$ that satisfies Constraint (2), the joint optimization problem reduces to the multi-price optimization problem under the PCL model, which is well studied in [16]. In the price optimization problem, [16] proves that the optimal prices are identical for all products. That is, letting $p_i^*(\mathbf{x})$ be the optimal price for product $i$, we have that $p_i^*(\mathbf{x}) = p^*(\mathbf{x})$ holds for any $i \in N$. Furthermore, the optimal expected revenue with given $\mathbf{x}$ can be characterized using $p^*(\mathbf{x})$, which leads to the following lemma.

**Lemma 1.** With given feasible assortment $\mathbf{x}$, the optimal price to maximize $R(p^*(\mathbf{x}), \mathbf{x})$ can be represented by

$$p^*(\mathbf{x}) = \frac{1 + W(A(\mathbf{x})/e)}{\beta}, \quad (4)$$

where $A(\mathbf{x})$ is calculated trough

$$A(\mathbf{x}) = \sum_{i=1}^{n-1} \sum_{j=i+1}^{n} (\exp(\alpha_i/\gamma_{ij})x_i + \exp(\alpha_j/\gamma_{ij})x_j)^{\gamma_{ij}}. \quad (5)$$

And the optimal expected revenue is $R(p^*(\mathbf{x}), \mathbf{x}) = W(A(\mathbf{x})/e)/\beta$, where $W(y)$ is the Lambert-W function to solve $we^w = y$ for $w$.

The proof of Lemma 1 is omitted here to save space since it is a direct generalization of Theorem 1 in [16] accounting for given assortment $\mathbf{x}$. Due to the fact that $R(p^*(\mathbf{x}), \mathbf{x}) = p^*(\mathbf{x}) - 1/\beta$, maximizing the expected revenue is equivalent to maximizing $p^*(\mathbf{x})$. So, for any given assortment $\mathbf{x}$, the optimal prices for all products are identical and can be calculated by (4), involving Lambert-W function. Since Lambert-W function is a strictly increasing function in $(0, +\infty)$, the optimal solution to maximize $A(\mathbf{x})$ subject to a capacity constraint must be the optimal assortment to the CAOPP. That is, we can further propose a two-step procedure to obtain the optimal solution to (1)-(3): In the first step, we find the feasible assortment $\mathbf{x}^*$ with the largest value of $A(\mathbf{x}^*)$ subject to the capacity limit; In the second step, we calculate $p^*(\mathbf{x}^*)$ based on Lemma 1. As a result, $(p^*(\mathbf{x}^*), \mathbf{x}^*)$ must be the optimal solution to the CAOPP. So, all that remains is to find $\mathbf{x}^*$, which can be obtained by solving the following non-linear integer program:

$$\max A(\mathbf{x}) = \sum_{i=1}^{n-1} \sum_{j=i+1}^{n} (\exp(\alpha_i/\gamma_{ij})x_i + \exp(\alpha_j/\gamma_{ij})x_j)^{\gamma_{ij}} \quad (6)$$

$$\text{s.t.} \quad \sum_{i=1}^{n} w_i x_i \leq C; \quad x_i \in \{0,1\}, \forall i \in N.$$

To solve the above problem, we first linearize the objective function by using the characteristic of the PCL model that there are exact two products under each nest. To do so, we follow the similar approach used in [17-18]. Letting $\rho_{ij} = (\exp(\alpha_i/\gamma_{ij}) + \exp(\alpha_j/\gamma_{ij}))^{\gamma_{ij}}$ and $\theta_i = \exp(\alpha_i)$, $A(\mathbf{x})$ can be re-written as $A(\mathbf{x}) = \sum_{i=1}^{n-1} \sum_{j=i+1}^{n} \rho_{ij} x_i x_j + \theta_i x_i (1 - x_j) + \theta_j (1 - x_i) x_j = \sum_{i=1}^{n-1} \sum_{j=i+1}^{n} \mu_{ij} x_i x_j + \theta_i x_i + \theta_j x_j$, where $\mu_{ij} = \rho_{ij} - \theta_i - \theta_j < 0$. Replacing the product $x_i x_j$ with $y_{ij} \geq 0$, we can ensure that $y_{ij} = 1$ if and only if $x_i = x_j = 1$ by introducing additional constraints $1 + y_{ij} \geq x_i + x_j, x_i \geq y_{ij}$, and $x_j \geq y_{ij}$. Using the fact that $\mu_{ij} < 0$, it is shown in [17] that the latter two constraints are redundant. Now, we can formulate our mixed integer program with linear objective function to solve (6) as follows:

$$\max A(\mathbf{x}) = \sum_{i=1}^{n-1} \sum_{j=i+1}^{n} \mu_{ij} y_{ij} + \theta_i x_i + \theta_j x_j \quad (7)$$

$$\text{s.t.} \quad \sum_{i=1}^{n} w_i x_i \leq C;$$
$$1 + y_{ij} \geq x_i + x_j, \forall i \neq j \in N;$$
$$y_{ij} \geq 0, x_i \in \{0,1\}, \forall i \in N.$$

The advantage of this formulation is that any MIP software package can be implemented to obtain the optimal assortment $\mathbf{x}^*$ with given capacity limit. However, solving an MIP to optimality can be quite time-consuming when there are thousands of products to consider since it is NP-hard. This complexity result prompts us to develop efficient heuristic algorithms to obtain high-quality solutions to (7) quickly.

We first use the linear relaxation of the MIP in (7), whose optimal objective function is denoted as $A(\overline{\mathbf{x}})$, to bound the value of $A(\mathbf{x}^*)$, which is used to gauge the quality of the solutions obtained by our heuristics. Since the optimal expected revenue is a non-decreasing with respect to $A(\mathbf{x})$, we can have $p^*(\overline{\mathbf{x}}) = (1 + W(A(\overline{\mathbf{x}})/e)))/\beta > p^*(\mathbf{x}^*)$ and $\overline{R} = p^*(\overline{\mathbf{x}}) - 1/\beta > R(p^*(\mathbf{x}^*), \mathbf{x}^*)$. That is, we can use the linear relaxation of the MIP in (7) and Lemma 1 to obtain an upper bound of CAOPP.

In our basic heuristic algorithm, a feasible solution is constructed as follows: (a) Sort the candidate products in decreasing order according to the values of $\left\{\frac{\exp(\alpha_i)}{w_i} : i \in N\right\}$; (b) Greedily add the feasible products, satisfying the capacity limits, into the current assortment until there is no feasible products that can be added into the assortment. As a result, we can obtain a feasible solution, denoted as $\mathbf{x}_H$, within $O(n\log n + n) = O(n\log n)$ time.

To further improve the quality of the obtained heuristic solution, we introduce a meta-heuristic algorithm based on a Greedy Randomized Adaptive Search Procedure, referred to as GRASP. In our GRASP, a feasible solution is randomized constructed based on our greedy heuristic and local perturbations are followed to get a local optimal solution. To do so, we introduce a parameter $\alpha \in [1, \text{Alpha}]$ in the second procedure (b) of our heuristic. That is, we randomly add the feasible product in the list composed of the $\alpha$ best indices according to the values of $\exp(\alpha_i)/w_i$ for any product $i$, which has not been selected in the current assortment, until no feasible products exist. Then, in the local search phase, we randomly choose two variables $x_i, x_j (x_i \neq x_j)$ and flip their values to zero or one. If the new solution is feasible and improves the objective function, we store the change and keep it as the current solution. Otherwise, we continue the random perturbation until the prescribed maximum number of iteration is reached. For ease of reading, the pseudo-code for our GRASP is shown in Procedure GRASP($Alpha, MaxIter$). As a result, we can obtain a feasible solution, denoted as $\mathbf{x}_{GRASP}$ and we can have $A(\mathbf{x}_{GRASP}) \geq A(\mathbf{x}_H)$ because when $\alpha = 1$, the construction phase is the same with the heuristic algorithm.

Procedure GRASP($Alpha, MaxIter$)
    For $\alpha = 1, 2, \cdots, Alpha$ do
        Obtain a feasible solution with heuristic algorithm
        With parameter $\alpha$
        For iteration=$1, 2, \cdots, MaxIter$ do
            Local search for better solutions by flipping
            the values of $x_i, x_j$
        End For
    End For

To sum up, we can now solve the MIP in (7) to obtain the optimal assortment, apply the greedy heuristic algorithm as well as the GRASP to obtain high-quality assortment with reasonable time consumed, and solve the linear relaxation of (7) to get an upper bound of $A(x)$. Then, for these four approaches, we can further use Lemma 1 to obtain the corresponding optimal price for all products as well as the corresponding expected revenue obtained from each customer. Therefore, we can have $R(p^*(\bar{\mathbf{x}}), \bar{\mathbf{x}}) \geq R(p^*(\mathbf{x}^*), \mathbf{x}^*) \geq R(p^*(\mathbf{x}_{GRASP}), \mathbf{x}_{GRASP}) \geq R(p^*(\mathbf{x}_H), \mathbf{x}_H)$ based on the discussions in this section. Now, we try to test the numerical performance of our approaches in the next section.

## IV. RESULTS

In this section, we compare the numerical performance of proposed algorithms in Section III based on the quality of the solutions obtained and the running time consumed. Throughout our numerical experiments, we use the following approach for generating problem instances. We vary the number of total products over $n \in \{400, 600, 800, 1000\}$. We set the capacity limit $C = \kappa \sum_{i=1}^{n} w_i$, where $\kappa$ is a fraction of the total space consumption of all products and varies over $\kappa \in \{0.02, 0.04, 0.06\}$. This setup provides 12 parameter combinations for $(n, \kappa)$. For each combination, we generate 25 random problem instance, each of which we sample the parameters as follows: we sample $\{\exp(\alpha_i): i \in N\}$ from the uniform distribution $U(0,5]$; the weight $w_i$ of product $i$ is generated from the uniform distribution $U[1,10]$; and the scale parameter $\gamma_{ij}$ for nest $\langle i, j \rangle$ is generated from the uniform distribution $U[0.1, 1]$ and $\gamma_{ij} = \gamma_{ji}$. Moreover, we set price sensitivity $\beta = 0.1$. The value of Alpha and the number of maximum iteration in local search used in our GRASP are set to 5 and 80, respectively. Then, for each problem instance, we obtain the optimal solution through solving MIP and Lambert-W function. We also apply the greedy heuristic algorithm and the GRASP to obtain high-quality solutions. Our numerical experiments are conducted on a Windows server with CPU of 2.00 GHz and RAM of 512 GB. The algorithms are all coded in Matlab linked to the CPLEX 12.6 optimization routines.

The numerical results for solving the CAOPP under the PCL model is summarized in Table I. Column 1 shows the parameter combination for instances in the form of $(n, \kappa)$. Let $Opt_{obj}^k, UB_{obj}^k, He_{obj}^k$, and $GRASP_{obj}^k$ be the expected revenue obtained by the MIP in (7), the linear relaxation of the MIP, the heuristic algorithm, and the GRASP, respectively, for problem instance $k \in \{1,2,\cdots,25\}$ with given parameter combination. Columns 2-3 report the average and maximum running time in 1 second to obtain the optimal solution based on solving the MIP and Lemma 1, respectively. Similarly, Columns 4-9, respectively, report the average and maximum running time in 1 second of obtaining an upper bound through linear relaxation of the MIP, a feasible solution through our greedy heuristic algorithm, and a solution through the GRASP. Columns 10-11 report the average and maximum optimality gap with MIP, defined as $\left(1 - \frac{Opt_{obj}^k}{UB_{obj}^k}\right) \times 100\%$, respectively. The average and maximum optimality gaps with heuristic algorithm and GRASP, defined as $\left(1 - \frac{He_{obj}^k}{UB_{obj}^k}\right) \times 100\%$ and $\left(1 - \frac{GRASP_{obj}^k}{UB_{obj}^k}\right) \times 100\%$, respectively, are reported in Columns 12-15 correspondingly.

Shown in Table I, the running time for our optimal approach based on MIP increases dramatically as the number of products and $\kappa$ increase. For combination (1000, 0.06), the maximum running time can reach 2263.0 seconds. On the other hand, the heuristic algorithm can obtain a feasible solution within 0.764 second. With a local search phase, the GRASP does need more time to terminate. Meanwhile, solving a linear relaxation for upper bounds seems to be a practical way to gauge the quality of solutions obtained by other algorithms with limited time consumed. As for the optimality gaps, we found that the average optimality gaps with MIP over all 300 instances is only 0.03% and the maximum is only 0.20%. That is, our approach to obtain an upper bound of the joint optimization problem is numerically tight. Meanwhile, the optimality gaps with our heuristic is 0.18% on average and 0.46% at the worst case. That is, our greedy heuristic algorithm can generate a considerable high-quality solution much more efficiently compared to the approach based on MIP. As for the GRASP, the numerical results indicate that it does improve the solution quality, especially, it can significantly improve the solution quality at the worst case. As we can see from Table 1, the optimality gap at the worst case is reduced from 0.46% to 0.31% compared with the naïve heuristic algorithm. To sum up, our heuristic algorithms are quite efficient to obtain high-quality feasible solutions with reasonable time consumed.

TABLE I
NUMERICAL RESULTS FOR DIFFERENT APPROACHES

| Para. comb. $(n, \kappa)$ | Run. time (s) for MIP | | Run. time (s) for UB | | Run. time (s) for Heuristic | | Run. time (s) for GRASP | | Gap (%) with MIP | | Gap (%) with Heuristic | | Gap (%) with GRASP | |
|---|---|---|---|---|---|---|---|---|---|---|---|---|---|---|
| | Avg. | Max. | Avg. | Max. | Avg. | Max. | Avg. | Max. | Avg. | Max. | Avg. | Max. | Avg. | Max. |
| (400, 0.02) | 9.9 | 11.8 | 3.6 | 5.6 | 0.150 | 0.176 | 39.1 | 39.7 | 0.11 | 0.20 | 0.28 | 0.46 | 0.20 | 0.31 |
| (400, 0.04) | 17.5 | 28.7 | 11.6 | 19.4 | 0.153 | 0.167 | 40.1 | 40.6 | 0.05 | 0.08 | 0.20 | 0.30 | 0.18 | 0.23 |
| (400, 0.06) | 38.5 | 50.2 | 32.3 | 51.7 | 0.156 | 0.173 | 40.9 | 41.2 | 0.03 | 0.05 | 0.22 | 0.26 | 0.21 | 0.25 |
| (600, 0.02) | 27.0 | 33.2 | 12.3 | 14.5 | 0.281 | 0.322 | 88.9 | 89.6 | 0.06 | 0.10 | 0.17 | 0.27 | 0.13 | 0.20 |
| (600, 0.04) | 220.7 | 228.8 | 57.7 | 84.4 | 0.285 | 0.303 | 91.2 | 91.9 | 0.03 | 0.05 | 0.17 | 0.23 | 0.16 | 0.19 |
| (600, 0.06) | 275.2 | 285.7 | 141.6 | 194.1 | 0.296 | 0.320 | 92.9 | 94.0 | 0.02 | 0.03 | 0.19 | 0.23 | 0.18 | 0.22 |
| (800, 0.02) | 84.5 | 691.3 | 30.1 | 35.6 | 0.458 | 0.505 | 161.7 | 163.5 | 0.05 | 0.08 | 0.15 | 0.22 | 0.12 | 0.18 |
| (800, 0.04) | 699.4 | 733.8 | 159.7 | 215.3 | 0.470 | 0.500 | 165.2 | 167.2 | 0.02 | 0.03 | 0.15 | 0.19 | 0.13 | 0.16 |
| (800, 0.06) | 827.2 | 856.5 | 395.2 | 480.0 | 0.483 | 0.535 | 168.1 | 169.9 | 0.01 | 0.02 | 0.17 | 0.20 | 0.16 | 0.19 |
| (1000, 0.02) | 244.9 | 1778.0 | 68.1 | 80.0 | 0.663 | 0.698 | 246.4 | 249.3 | 0.04 | 0.06 | 0.13 | 0.19 | 0.11 | 0.13 |
| (1000, 0.04) | 1919.3 | 2134.3 | 379.0 | 475.5 | 0.674 | 0.732 | 250.4 | 256.1 | 0.01 | 0.02 | 0.14 | 0.17 | 0.13 | 0.15 |
| (1000, 0.06) | 2143.5 | 2263.0 | 951.5 | 1166.7 | 0.717 | 0.764 | 257.6 | 263.6 | 0.01 | 0.02 | 0.16 | 0.19 | 0.16 | 0.18 |
| Average | ---- | ---- | ---- | ---- | ---- | ---- | ---- | ---- | 0.03 | 0.06 | 0.18 | 0.24 | 0.16 | 0.20 |

## V. CONCLUSION

In this paper, we investigate the capacitated assortment optimization problem with pricing. The customers' purchase behavior is modeled using the paired combinatorial logit model. We first formulate this problem as a non-linear mixed integer program and propose a two-step approach to obtain the optimal solution based on solving an MIP and a Lambert-W function. To further improve the numerical performance, we develop a greedy heuristic and a GRASP to obtain high-quality solutions quickly. The numerical experiments indicate that our heuristic can obtain high-quality solutions whose deviations from the upper bound are only 0.16% on average and 0.48% at the worst case. We now outline two possible future research directions: (a) In our settings, the price sensitivity is identical for all products. However, there can be scenarios where customers have different sensitivities over different products. It would be of interest to investigate the case with different price sensitivities; (b) When the retailer can determine the prices and the assortment dynamically over finite selling periods, our algorithms may not remain efficient. Therefore, it is of great interest to study the corresponding problem in a dynamic environment.

## REFERENCES

[1] K. D. Chen, W. H. Hausman, "Technical note: Mathematical properties of the optimal product line selection problem using choice-based conjoint analysis", *Management Science*, vol. 46, no. 2, pp. 327–332, 2000.
[2] R. X. Wang, "Capacitated assortment and price optimization under the multinomial logit model", *Operations Research Letters*, vol. 40, no. 6, pp. 492–497, 2012.
[3] G. Gallego, H. Topaloglu, "Constrained assortment optimization for the nested logit model", *Management Science*, vol. 60, no. 10, pp. 2583–2601, 2014.
[4] J. Davis, G. Gallego, H. Topaloglu, "Assortment planning under the multinomial logit model with totally unimodular constraint structures", unpublished.
[5] O. Besbes, D. Saure, "Product assortment and price competition under multinomial logit demand", *Production and Operations Management*, vol. 25, no. 1, pp. 114–127, 2016.
[6] K. E. Train, "Discrete choice methods with simulation", *Cambridge university press*, 2009.
[7] H. C. Williams, "On the formation of travel demand models and economic evaluation measures of user benefit", *Environment and planning A*, vol. 9, no. 3, pp. 285–344, 1977.
[8] A. G. Kok, Y. Xu, "Optimal and competitive assortments with endogenous pricing under hierarchical consumer choice models", *Management Science*, vol. 57, no. 9, pp. 1546–1563, 2011.
[9] W. Z. Rayfield, P. Rusmevichientong, H. Topaloglu, "Approximation methods for pricing problems under the nested logit model with price bounds", *Informs Journal on Computing*, vol. 27, no. 2, pp. 335–357, 2015.
[10] R. Chen, H. Jiang, "Capacitated assortment and price optimization for customers with disjoint consideration sets", *Operations Research Letters*, vol. 45, no. 2, pp. 170–174, 2017.
[11] S. Jagabathula, P. Rusmevichientong, "A nonparametric joint assortment and price choice model", *Management Science*, vol. 63, no. 9, pp. 3128–3145, 2016.
[12] F. S. Koppelman, C.-H. Wen, "The paired combinatorial logit model: Properties, estimation and application", *Transportation Research Part B: Methodological*, vol. 34, no. 2, pp. 75–89, 2000.
[13] J. Prashker, S. Bekhor, "Investigation of stochastic network loading procedures", *Transportation Research Record: Journal of the Transportation Research Board*, no. 1645, pp. 94–102, 1998.
[14] A. Chen, S. Ryu, X. Xu, K. Choi, "Computation and application of the paired combinatorial logit stochastic user equilibrium problem", *Computers & Operations Research*, vol. 43, pp. 68–77, 2014.
[15] A. Karoonsoontawong, D.-Y. Lin, "Combined gravity model trip distribution and paired combinatorial logit stochastic user equilibrium problem", *Networks and Spatial Economics*, vol. 15, no. 4, pp. 1011–1048, 2015.
[16] H. Li, S. Webster, "Optimal pricing of correlated product options under the paired combinatorial logit model", *Operations Research*, vol. 65, no. 5, pp. 1215–1230, 2017.
[17] H. Zhang, P. Rusmevichientong, H. Topaloglu, "Assortment optimization under the paired combinatorial logit model", unpublished.
[18] J. Feldman, "Space constrained assortment optimization under the paired combinatorial logit model", unpublished.